\documentclass{amsart}%
\usepackage{amsfonts}
\usepackage[hiresbb]{graphicx}
\usepackage{amsmath}
\usepackage{amssymb}
\usepackage{version}
\usepackage{epsfig}
\usepackage{epstopdf}
\usepackage{bibmods}%
\providecommand{\U}[1]{\protect\rule{.1in}{.1in}}
\DeclareGraphicsExtensions{.eps}
\newtheorem{defn}{Definition}[section]
\newtheorem{thm}[defn]{Theorem}
\newtheorem{exmp}[defn]{Example}

\newtheorem{rmk}[defn]{Remark}
\newtheorem{prop}[defn]{Proposition}
\newtheorem{cor}[defn]{Corollary}

\begin{document}
\title{Symbolic Dynamics of Odd Discontinuous Bimodal Maps}
\author{Henrique M. Oliveira}
\address{Center for Mathematical Analysis, Geometry and Dynamical Systems, Mathematics
Department, Instituto Superior T\'{e}cnico, Universidade de Lisboa, Av.
Rovisco Pais, 1049-001 Lisboa, Portugal }
\email{holiv@math.ist.ulisboa.pt}

\begin{abstract}
Iterations of odd piecewise continuous maps with two discontinuities, i.e.,
symmetric discontinuous bimodal maps, are studied. Symbolic dynamics is
introduced. The tools of kneading theory are used to study the homology of the
discrete dynamical systems generated by the iterations of that type of maps.
When there is a Markov matrix, the spectral radius of this matrix is the
inverse of the least root of the kneading determinant.

\end{abstract}
\maketitle

\section{Introduction}

In this paper we apply techniques of Markov partitions and kneading theory to
the study of iterates of discontinuous maps of the interval (or the real line)
in itself. We show that these systems can be studied with a proper framework,
which is related to kneading theory and Markov matrices.

We cite, as examples of discontinuous one dimensional cases, the Lorenz maps,
Newton maps, circle and tree maps, see \cite{Al,Alv,CUR,Glen} among other literature.

In \cite{LSR} Lampreia and Sousa Ramos studied symbolic dynamics of continuous
bimodal maps on the compact interval. Using similar techniques, we study in
this paper the case of symmetric (odd) discontinuous maps in the real line or
some suitable interval with two discontinuity points and three maximal
intervals (laps) of continuity, which are as well maximal intervals of
monotonicity. We call to this type of mapping a \textit{symmetric bimodal
discontinuous map} because of the existence of exactly three laps as in the
continuous bimodal case.

In section two, we introduce the notation, the main definitions and revision
of basic results. We include as well, the notions of symbolic dynamics,
kneading theory and Markov partitions. We relate these concepts with lap
growth number. We tried to define with great detail all the concepts
presented. Since good definitions are essential for the constructive proof of
the main result, which is actually done along the full length of the paper,
section two is relatively long.

In section three, we present the main result of the paper, i.e., the spectral
radius of the Markov matrix is the inverse of the least root of the kneading
determinant for that kind of maps. We point out that the introduction of the
linear operator $\mu$ in section three is one of the main ideas of this paper
along with the matrix $\Theta$ relating the kneading and the Markov data. The
linear transformation $\mu$, representing the symmetry of this type of
non-continuous maps, is completely different of its continuous counterpart
\cite{LSR}. We think that the proof of the result can be instructive giving
methods that can be applied to other non-continuous mappings.

\subsection{Motivation}

The iterates of the complex tangent family $\lambda\tan z$, introduced in
\cite{dev-keen} and \cite{keen-kotus},\ when the parameter $\lambda=i\beta$ is
pure imaginary and the initial condition $x_{0}$ is a real number can be
identified with the iterates of the real alternating map $\left[  f_{1,\beta
},f_{2,\beta}\right]  $ \cite{EH} in the real line%
\begin{align*}
x_{1}  &  =f_{1,\beta}\left(  x_{0}\right)  =\beta\tan\left(  x_{0}\right) \\
x_{2}  &  =f_{2,\beta}\left(  x_{1}\right)  =-\beta\tanh\left(  x_{1}\right)
\\
x_{3}  &  =f_{1,\beta}\left(  x_{2}\right) \\
x_{4}  &  =f_{2,\beta}\left(  x_{3}\right)  .\\
&  \vdots
\end{align*}
The composition map $g_{\beta}$ is%
\begin{equation}%
\begin{array}
[c]{cccc}%
g_{\beta}: & \mathbb{R} & \mathbb{\mapsto} & \mathbb{R}\\
& x & \rightarrow & -\beta\tanh\left(  \beta\tan\left(  x\right)  \right)  ,
\end{array}
\label{comp}%
\end{equation}
which can be interpreted as the second return map to the real axis for the
mapping $\lambda\tan\left(  z\right)  $. Knowing $x_{0}$ and $g_{\beta}$, we
obtain all the even iterates of the original system. To obtain the odd
iterates knowing the even iterates is easy%
\[
x_{2n+1}=\beta\tan\left(  x_{2n}\right)  .
\]

The geometric behavior of the maps $g_{\beta}$ in this family depends on the
parameter $\beta$. The map is periodic and the real line is mapped on the
interval $I=\left(  -\beta,\beta\right)  $. We restrict the map only to the
interval $I$. When $\frac{\pi}{2}<\beta<\frac{3\pi}{2}$ the maps $g_{\beta}$
have two discontinuities. The study of the real projection of the complex
tangent map is a good clue to the dynamics in the complex plane, similarly to
the case of quadratic maps.

In this paper, we center our study on the symbolic dynamics of the iterates of
maps $F$ with the same geometrical properties of $g_{\beta}$. Considering that
the tangent family was an initial motivation and a good example, we point out
that the results are independent on the choice of the family.

\section{Basics}

\subsection{Bimodal symmetric discontinuous map}

\begin{defn}
\label{DEF1} Bimodal symmetric discontinuous map of type $\left(
-,-,-\right)  $. Let $I=\left(  -a,a\right)  $ be a real interval (where $a$
can be $+\infty$) and $F:I\mathbb{\mapsto}I$, such that:

\begin{enumerate}
\item \label{Simetrica}$F$ is odd $F\left(  x\right)  =-F\left(  -x\right)  $

\item $F$ is piecewise continuous having two discontinuities $c_{1}<c_{2}$,
$c_{1}=-c_{2}$, where $\lim_{x\rightarrow c_{i}^{\pm}}F\left(  x\right)  =\pm
a$,\ and$\ \lim_{x\rightarrow\pm a}F\left(  x\right)  =\pm b$, where $b$ is a
real number.

\item $F$ is decreasing in every interval of continuity $\left(
-a,c_{1}\right)  $, $\left(  c_{1},c_{2}\right)  $ and $\left(  c_{2}%
,a\right)  $.
\end{enumerate}
\end{defn}

\begin{exmp}
The family of maps $g_{\beta}$ defined in (\ref{comp}) is a family of bimodal
symmetric discontinuous maps.
\end{exmp}

Definition \ref{DEF1} applies to maps with infinite jumps at $c_{1}$ and
$c_{2}$ as we can see in the next example. Actually, as we see in the next
example, any such map is smoothly conjugated to a map with finite jumps via a diffeomorphism.

\begin{exmp}
\label{Galpha}Consider $u:\mathbb{R\mapsto R}$ such that
\[
u\left(  x\right)  =\left\{
\begin{array}
[c]{rrr}%
-1 & \text{if} & x\leq-\frac{1}{2}\\
0 & \text{if} & -\frac{1}{2}<x<\frac{1}{2}\\
1 & \text{if} & \frac{1}{2}\leq x
\end{array}
\right.  ,
\]
the family
\[%
\begin{array}
[c]{cccc}%
G_{\alpha}: & \mathbb{R} & \rightarrow & \mathbb{R}\\
& x & \mathbb{\mapsto} & \frac{x}{4x^{2}-1}-\alpha\text{ }u\left(  x\right)
\end{array}
\]
is a family of bimodal symmetric discontinuous maps with $b=-\alpha$,
$c_{1}=-\frac{1}{2}$, $c_{2}=\frac{1}{2}$ and $a=+\infty$. Any map in this
family satisfies definition \ref{DEF1} and is smoothly conjugated to a map
with finite jumps using for instance the diffeomorphism $h\left(  x\right)
=\arctan\left(  x\right)  $ such that%
\[
\widetilde{G}_{\alpha}\left(  x\right)  =h\circ G_{\alpha}\circ h^{-1}\left(
x\right)  \text{, }x\in\left(  -\frac{\pi}{2},\frac{\pi}{2}\right)  \text{.}%
\]
The map $\widetilde{G}_{\alpha}$ can be prolonged by continuity to the
endpoints $\pm\frac{\pi}{2}$ of the interval, since
\[
\lim_{x\rightarrow\pm\frac{\pi}{2}}h\circ G_{\alpha}\circ h^{-1}\left(
x\right)  =\lim_{x\rightarrow\pm\infty}h\circ G_{\alpha}\left(  x\right)
=\mp\arctan\left(  \alpha\right)  .
\]

\end{exmp}

\subsection{Symbolic dynamics}

For sake of completeness and readability we introduce here briefly notions
well known like orbit, periodic orbit and symbolic itinerary among other
concepts, see for instance \cite{DS}.

\begin{defn}
We define the orbit of a real point $x_{0}$ as a sequence of numbers $O\left(
x_{0}\right)  =\left\{  x_{j}\right\}  _{j=0,1,...}$ such that $x_{j}%
=F^{j}\left(  x_{0}\right)  $ where $F^{j}$ is the $j$-th composition of $F$
with itself.
\end{defn}

\begin{defn}
Any point $x$ is periodic with period $n>0$ if the condition $F^{n}\left(
x\right)  =x$ is fulfilled with $n$ minimal.
\end{defn}

Because of condition \ref{Simetrica} in definition \ref{DEF1},\ the orbit of
any point $x$ is symmetric relative to the orbit of $-x$. To avoid ambiguities
in the definition of the orbit of the pre discontinuity points we adopt the
convention that $F\left(  c_{1}\right)  =F\left(  c_{1}^{-}\right)  =-a$ and
$F\left(  c_{2}\right)  =F\left(  c_{2}^{+}\right)  =+a$.

\begin{defn}
\cite{OLI} Consider the alphabet $\mathcal{A}=\left\{  L,A,M,B,R\right\}  $
the address $\mathtt{A}\left(  x\right)  $ of a real point $x$ is defined such
that
\[
\mathtt{A}\left(  x\right)  =\left\{
\begin{array}
[c]{ccc}%
L & \text{if} & x<c_{1}\\
A & \text{if} & x=c_{1}\\
M & \text{if} & c_{1}<x<c_{2}\\
B & \text{if} & c_{2}=x\\
R & \text{if} & c_{2}<x
\end{array}
\right.  .
\]

\end{defn}

We can apply this function to an orbit of a given real point $x_{0}$, we
associate to that orbit one infinite symbolic sequence.

\begin{defn}
Consider the sequence of symbols in $\mathcal{A}$
\[
\operatorname{It}\left(  x_{0}\right)  =\mathtt{A}\left(  x_{0}\right)
\mathtt{A}\left(  x_{1}\right)  \mathtt{A}\left(  x_{2}\right)  ...\mathtt{A}%
\left(  x_{n}\right)  ...
\]
this infinite sequence is the symbolic itinerary of $x_{0}$.
\end{defn}

The orbit $O\left(  -a\right)  $ is
\[
\left\{  x_{j}^{\left(  1\right)  }:x_{j}^{\left(  1\right)  }=F^{j}\left(
-a\right)  ,j=0,1,\ldots\right\}  .
\]
The orbit of $O\left(  +a\right)  $ is
\[
\left\{  x_{j}^{\left(  2\right)  }:x_{j}^{\left(  2\right)  }=F^{j}\left(
+a\right)  ,j=0,1,\ldots\right\}  ,
\]
with $F\left(  +a\right)  =b$.

\begin{defn}
\cite{OLI} \emph{Kneading sequences and kneading pairs.} The kneading
sequences are defined as the symbolic itineraries of the orbits of $a$ and
$-a$. The kneading pair is the ordered pair formed by these two symbolic
sequences%
\[
\left(  \operatorname{It}\left(  a\right)  ,\operatorname{It}\left(
-a\right)  \right)  .
\]

\end{defn}

\begin{defn}
\emph{Order relation in }$\mathcal{A}$\emph{. }The order on $\mathcal{A}$\ is
naturally induced from the order in the real axis
\[
L\prec A\prec M\prec B\prec R.
\]

\end{defn}

\begin{defn}
\cite{OLI} \emph{Parity function} $\rho\left(  S\right)  $. Given any finite
sequence $S$ with length $p$, $\rho\left(  S\right)  $ is such that%
\[
\rho\left(  S\right)  =\left(  -1\right)  ^{p}.
\]

\end{defn}

\begin{defn}
\cite{OLI} Let $\mathcal{A}^{\mathbb{N}}$ denote the set of all sequences
written with the alphabet $\mathcal{A}$. We define an ordering $\prec$ on the
set $\mathcal{A}^{\mathbb{N}}$ such that: given two symbolic sequences
$P=P_{0}P_{1}P_{2}...$ and $Q=Q_{0}Q_{1}Q_{2}...$ let $n$ be the first integer
such that $P_{n}\neq Q_{n}$. Denote by $S=S_{0}S_{1}S_{2}...S_{n-1}$ the
common first subsequence of both $P$ and $Q$. Then, we say that $P\prec Q$ if
$P_{n}\prec Q_{n}$ and $\rho\left(  S\right)  =+1$ or $Q_{n}\prec R_{n}$ if
$\rho\left(  S\right)  =-1$. If no such $n$ exists then $P=Q$.
\end{defn}

This ordering is originated by the fact that when $x<y$ then
$\operatorname{It}\left(  x\right)  \preceq\operatorname{It}\left(  y\right)
$.

To state the rules of admissibility the shift operator $\sigma$ will be used,
defined as usual.

\begin{defn}
\emph{Shift operator }$\sigma$. The shift operator is defined%
\[
\sigma\left(  P_{0}P_{1}P_{2}...\right)  =P_{1}P_{2}...\text{.}%
\]
When we have a finite sequence $S$ the shift operator acts such that%
\[
\sigma\left(  S_{0}S_{1}S_{2}...S_{n-1}\right)  =S_{1}S_{2}...S_{n-1}%
S_{0}\text{.}%
\]

\end{defn}

The orbit of $+a$ has the symbolic itinerary $\operatorname{It}\left(
+a\right)  $. The sequence $\operatorname{It}\left(  +a\right)  $ is maximal
(resp. $\operatorname{It}\left(  -a\right)  $ is minimal) in the ordering
defined in this section. Maximal in the sense that every shift of the sequence
$\operatorname{It}\left(  +a\right)  $ is less or equal than
$\operatorname{It}\left(  +a\right)  $. Every orbit with initial condition
$x_{0}$ is symmetric to the orbit with initial condition $-x_{0}$. Thus any
orbit beginning by $+a$ is accompanied by a symmetric orbit started by $-a$.
Therefore, we shall focus the admissibility rules for kneading sequences only
on the itineraries with the first symbol $R$ (corresponding to $+a$).

\begin{defn}
\cite{OLI} \emph{Operator }$\tau.$ The operator $\tau:\mathcal{A}^{\mathbb{N}%
}\mapsto\mathcal{A}^{\mathbb{N}}$ is defined such that
\[
\tau L=R,\tau A=B,\tau M=M,\tau B=A,\tau R=L\text{.}%
\]
Given a sequence $Q=Q_{0}Q_{1}Q_{2}...,$ $\tau$ acts such that%
\[
\tau Q=\tau Q_{0}\tau Q_{1}\tau Q_{2}...
\]

\end{defn}

The operator $\tau$ interchanges the symbols $L$ and $R,$\ letting the symbols
$M$ unchanged. For instance $\tau\left(  \left(  RLMR\right)  ^{\infty
}\right)  =\left(  LRML\right)  ^{\infty}$.

\begin{prop}
\cite{OLI} $\operatorname{It}\left(  x_{0}\right)  =\tau\operatorname{It}%
\left(  -x_{0}\right)  $.
\end{prop}

\smallskip\noindent\textbf{Proof.} Is a direct consequence of condition
\ref{Simetrica} in definition \ref{DEF1}. $\square$

Given any itinerary of $+a$ denoted by $S$, the corresponding itinerary of
$-a$ is $\tau S$. The kneading pair is $\left(  S,\tau S\right)  $. To know
the kneading sequence $S$, corresponding to the orbit of $+a$, is to know the
kneading pair. By some abuse of notation, sometimes (mainly in the examples)
we use only the kneading sequence $S$ instead of the kneading pair.

\begin{defn}
\cite{OLI} \emph{Admissibility rules: }Let $S$ be a given sequence of symbols
and $\left(  S,\tau S\right)  $ be a pair of sequences. $\left(  S,\tau
S\right)  $ is a kneading pair and $S$ is a kneading sequence, if $S$
satisfies the admissibility condition: $\tau S\preceq\sigma^{k}S\preceq S$,
for every integer $k$. The set of the admissible sequences is denoted by
$\Sigma\subset\mathcal{A}^{\mathbb{N}}$.
\end{defn}

\begin{defn}
Given a finite sequence $P$ with length $p$, the sequence $S=P^{\infty}$ is
called a $p$-periodic sequence.
\end{defn}

We will work sometimes only with $P$ instead of $P^{\infty}$ when there is no
danger of confusion.

\begin{defn}
\cite{OLI} A bistable periodic orbit contains both the orbit of $+a$ and the
orbit of $-a$. Any bistable orbit has an itinerary $S=P^{\infty}=\left(  Q\tau
Q\right)  ^{\infty}$ or shortly $P=Q\tau Q$
\end{defn}

As a consequence of the previous definition bistable orbits and associated
symbolic itineraries must have even period.

\subsection{Kneading theory}

In \cite{knead} were introduced the concepts of invariant coordinate, kneading
increments, kneading matrix and kneading determinant. We will use the
definitions of the cited work with the convenient adaptations for the
discontinuous case. We present here a brief exposition of the results obtained
applying kneading theory to this type of maps.

\begin{defn}
\emph{Invariant coordinate of an initial condition }$\theta_{x_{0}}\left(
t\right)  $.$\ $Is defined using the sequence $X=X_{0}X_{1}X_{2}%
\ldots=\operatorname{It}\left(  x_{0}\right)  $. Is the formal power series
\[
\theta_{x_{0}}\left(  t\right)  =\sum\limits_{k=0}^{+\infty}\left(  -1\right)
^{k}X_{k}t^{k}.
\]

\end{defn}

With the notation $\theta_{c_{i}^{\pm}}\left(  t\right)  =\lim_{x\rightarrow
c_{i}^{\pm}}\theta_{x}\left(  t\right)  $, for each discontinuity point, the
kneading increment is defined.

\begin{defn}
\emph{Kneading increment and kneading matrix.} The kneading increment is%
\[
\nu_{i}\left(  t\right)  =\theta_{c_{i}^{+}}\left(  t\right)  -\theta
_{c_{i}^{-}}\left(  t\right)  .
\]
This quantity is a formal power series measuring the discontinuity. After
collecting the terms associated to each symbol, and remarking that, in this
case, $c_{1}^{-}$ corresponds to $L$, $c_{1}^{+}$ corresponds to $M$,
$c_{2}^{-}$ corresponds to $M$ and $c_{2}^{+}$ corresponds to $R$, the
decomposition
\[
\nu_{i}\left(  t\right)  =N_{i1}\left(  t\right)  L+N_{i2}\left(  t\right)
M+N_{i3}\left(  t\right)  R
\]
is obtained. The kneading matrix is%
\[
N=\left[
\begin{array}
[c]{rrr}%
N_{11}\left(  t\right)  & N_{12}\left(  t\right)  & N_{13}\left(  t\right) \\
N_{21}\left(  t\right)  & N_{22}\left(  t\right)  & N_{23}\left(  t\right)
\end{array}
\right]  \text{.}%
\]

\end{defn}

\begin{defn}
\emph{Kneading determinant.} Omitting the $j$-th column of the kneading matrix
we compute the determinant $D_{j}$. The kneading determinant is
\[
D\left(  t\right)  =\frac{\left(  -1\right)  ^{j+1}D_{j}}{\left(  1+t\right)
}.
\]

\end{defn}

The denominator in the kneading determinant results from the fact that $F$ is
decreasing in the three intervals where this map is defined \cite{knead}. Note
that $D_{1}=-D_{2}=D_{3}$.

\begin{defn}
\label{FI}\cite{OLI} Given a sequence $X=X_{1}X_{2}\ldots$ we define a
function $\Phi:\mathcal{A}\longmapsto\left\{  -1,0,1\right\}  $, such that
\[
\Phi\left(  X_{i}\right)  =\left\{
\begin{array}
[c]{r}%
-1\\
0\\
1
\end{array}%
\begin{array}
[c]{llll}%
\text{if} & X_{i} & = & L,A\\
\text{if} & X_{i} & = & M\\
\text{if} & X_{i} & = & R,B
\end{array}
\right.  .
\]

\end{defn}

\begin{defn}
\cite{OLI} Given a sequence $X=X_{1}X_{2}\ldots$ from $\mathcal{A}%
^{\mathbb{N}}$\ we define a formal power series $u\left(  t\right)  $, such
that $u\left(  t\right)  =\sum\limits_{k=1}^{+\infty}\left(  -1\right)
^{k}\Phi\left(  X_{k}\right)  t^{k}$. When $X$ is finite with length $p\ $we
define the formal polynomial $u_{p}\left(  t\right)  =\sum\limits_{k=1}%
^{p}\left(  -1\right)  ^{k}\Phi\left(  X_{k}\right)  t^{k}$.
\end{defn}

Let $S=S_{1}S_{2}\ldots\in\Sigma$ be a kneading sequence with the kneading
pair $\left(  S,\tau S\right)  $, then the kneading determinant is given by.%
\begin{equation}
D\left(  t\right)  =\frac{1+2u\left(  t\right)  }{t+1}. \label{Kneading}%
\end{equation}

When $S=P^{\infty}$ is $p$-periodic, the expression of the kneading
determinant simplifies into $D\left(  t\right)  \left(  t+1\right)
=1+\frac{2u_{p}\left(  t\right)  }{1-\left(  -1\right)  ^{p}t^{p}}.$ When the
kneading sequence $S$ is bistable: $S=P^{\infty}$ with $P=Q\tau Q$ of period
$p=2q$, with the associated kneading pair $\left(  S,\tau S\right)  =\left(
Q\tau Q,\left(  \tau Q\right)  Q\right)  $ the kneading determinant is
$\left(  1+\frac{2u_{q}\left(  t\right)  }{1+\left(  -1\right)  ^{q}t^{q}%
}\right)  \frac{1}{1+t}$.

\subsection{Growth number}

The kneading determinant is essential in the computation of the growth number
of laps.

\begin{defn}
Lap number $\ell\left(  F^{n}\right)  $ is the number of maximal intervals of
continuity of each composition of $F$ with itself.
\end{defn}

\begin{defn}
The growth number is defined
\begin{equation}
\rho=\lim_{n\rightarrow\infty}\sqrt[n]{\ell\left(  F^{n}\right)  }.
\label{Mis}%
\end{equation}
\qquad\qquad
\end{defn}

\begin{rmk}
\cite{OLI} The growth number of $F$ can be computed using the relation
\[
\rho=\frac{1}{t_{0}},
\]
where $t_{0}$ is the least root in the unit interval of the kneading
determinant $D\left(  t\right)  $. The proof is provided defining the power
series $\Lambda\left(  t\right)  =\sum_{n\geq1}\ell\left(  F^{n}\right)
t^{n-1}$, where each coefficient is the lap number of the iterate $F^{n}$.
This new power series is closely related to the kneading determinant because
of the relation $\ \Lambda\left(  t\right)  =\frac{1}{t\left(  1-t^{2}\right)
D\left(  t\right)  }-\frac{1}{t}$.
\end{rmk}

\begin{exmp}
\label{ex3}The kneading sequence $\left(  RMR\right)  ^{\infty}$ corresponds
to the kneading determinant $D\left(  t\right)  =\frac{1-2t-t^{3}}{\left(
t+1\right)  \left(  1+t^{3}\right)  }$, which is realized for instance by
$g_{\beta}\left(  x\right)  $ with $\beta$ approximately $3.1588$ or by
$G_{\alpha}$ (from example \ref{Galpha}) with $\alpha=\frac{1}{4}\left(
\sqrt{5}-1\right)  $. We obtain $\Lambda\left(  t\right)  =3+7t+17t^{2}%
+39t^{3}+87t^{4}+193t^{5}+\ldots$, in this case $\ell\left(  F\right)  =3$,
$\ell\left(  F^{2}\right)  =7$, $\ell\left(  F^{3}\right)  =17\ldots$
\end{exmp}

It will be an interesting work to see if the usual relationship between the
topological entropy and growth number still remains valid in the case of
discontinuous maps.

\subsection{Markov partition}

Whenever we can define Markov matrices, the method of Markov transition
matrices in the case of continuous maps is an equivalent approach to
the\ computation the roots of the kneading determinant. To each $p$-periodic
kneading pair we associate a Markov transition matrix, see \cite{LSR} and
related references on that paper. Now, denote by
\begin{align*}
x_{j}^{\left(  2\right)  }  &  =F^{j}\left(  c_{2}^{+}\right)  ,\text{
}j=0,1,\ldots,p-1\text{,}\\
x_{j}^{\left(  1\right)  }  &  =F^{j}\left(  c_{1}^{-}\right)  ,\text{
}j=0,1,\ldots,p-1\text{,}%
\end{align*}
the orbits of the discontinuity points. An ordered sequence $\left(
z_{k}\right)  _{k=1,\ldots,2p}$ is obtained reordering the elements
$x_{j}^{\left(  m\right)  }$, $m=1,2$, and getting a partition
\[
I_{k}=\left(  z_{k},z_{k+1}\right)  \text{ with }k=1,\ldots,2p-1\text{.}%
\]
The discontinuity points are present in the above partition. We call
$z_{k_{1}}=c_{1}$ and $z_{k_{2}}=c_{2}$. To compute the Markov matrix note
that $I_{k_{1}-1}=\left(  z_{k_{1}-1},c_{1}^{-}\right)  $ and $I_{k_{1}%
}=\left(  c_{1}^{+},z_{k_{1+1}}\right)  $ and similarly with the two intervals
adjacent to the discontinuity point $c_{2}$. With this precision made, the
Markov transition matrix can be defined.

\begin{defn}
\emph{The Markov transition matrix }$\Psi=\left[  \psi_{ij}\right]  $ is
defined by the rule:
\[
\psi_{ij}=\left\{
\begin{array}
[c]{cc}%
1 & \text{if }I_{j}\subset F\left(  I_{i}\right)  \text{,}\\
0 & \text{otherwise.}%
\end{array}
\right.
\]

\end{defn}

In \cite{LSR} the relationship between Markov partitions and kneading theory
is explained for bimodal continuous maps. It is also presented the proof of
the equality of the reciprocal of $t_{0}$ and the spectral radius of the
matrix $\Psi$. In this paper we will prove the same equivalence of definitions
in the case of bimodal symmetric discontinuous maps.

In the next example we obtain this equivalence for a particular case,
computing directly both the kneading determinant and the characteristic
polynomial of the Markov matrix.

\begin{exmp}
The kneading pair
\[
\left(  \left(  RMR\right)  ^{\infty},\left(  LML\right)  ^{\infty}\right)
\]
corresponds to a pair of orbits satisfying
\begin{equation}
x_{1}^{\left(  1\right)  }<x_{0}^{\left(  1\right)  }=c_{1}<x_{2}^{\left(
2\right)  }<x_{2}^{\left(  1\right)  }<x_{0}^{\left(  2\right)  }=c_{2}%
<x_{1}^{\left(  2\right)  }, \label{order}%
\end{equation}
renaming the elements of the partition we get
\begin{align*}
z_{1}  &  =x_{1}^{\left(  1\right)  },z_{2}=x_{0}^{\left(  1\right)  }%
,z_{3}=x_{2}^{\left(  2\right)  },\\
z_{4}  &  =x_{2}^{\left(  1\right)  },z_{5}=x_{0}^{\left(  2\right)  }%
,z_{6}=x_{1}^{\left(  2\right)  }.
\end{align*}
The Markov matrix is
\[
\Psi=\left[
\begin{array}
[c]{rrrrr}%
1 & 1 & 1 & 0 & 0\\
0 & 0 & 0 & 0 & 1\\
0 & 1 & 1 & 1 & 0\\
1 & 0 & 0 & 0 & 0\\
0 & 0 & 1 & 1 & 1
\end{array}
\right]  .
\]
The smallest solution $t_{0}$ of the equation
\[
\det\left(  I-\Psi t\right)  =\left(  1-t+t^{2}\right)  \left(  1-2t-t^{3}%
\right)  =0
\]
in the unit interval gives the growth number $\rho=\frac{1}{t_{0}}=\frac
{1}{0.4534}=2.2056$, exactly the same root obtained with the kneading
determinant of the example \ref{ex3}.
\end{exmp}

\subsection{Relation between Markov partition and the orbits of the
discontinuity points}

Giving the $p$-periodic orbits
\[
O\left(  c_{2}^{+}\right)  =\left\{  x_{j}^{\left(  1\right)  }:x_{j}^{\left(
1\right)  }=F^{j}\left(  c_{2}^{+}\right)  ,\text{ }j=0,1,\ldots,p-1\right\}
\]
and
\[
O\left(  c_{1}^{-}\right)  =\left\{  x_{j}^{\left(  2\right)  }:x_{j}^{\left(
2\right)  }=F^{j}\left(  c_{1}^{-}\right)  ,\text{ }j=0,1,\ldots,p-1\right\}
,
\]
we define the vector
\[
y=\left[
\begin{array}
[c]{l}%
y_{1}\\
\vdots\\
y_{p}\\
y_{p+1}\\
\vdots\\
y_{2p}%
\end{array}
\right]  =\left[
\begin{array}
[c]{l}%
x_{0}^{\left(  2\right)  }\\
\vdots\\
x_{p-1}^{\left(  2\right)  }\\
x_{0}^{\left(  1\right)  }\\
\vdots\\
x_{p-1}^{\left(  1\right)  }%
\end{array}
\right]  .
\]
Let $z$ be the vector $\left\{  z_{i}\right\}  _{i=1,\ldots,2p}$ where
$z_{i-1}<z_{i}<z_{i+1}$ are the ordered elements of $y$. There is a
$2p\times2p$ permutation matrix $\pi$ such that
\[
z=\pi y.
\]
Let $\mathtt{x}_{k}^{\left(  2\right)  }=\operatorname{It}\left(
x_{k}^{\left(  2\right)  }\right)  $, the symbolic itinerary of $x_{k}%
^{\left(  2\right)  }$, for $k=0,\ldots,p-1$. It is clear that $\mathtt{x}%
_{1}^{\left(  2\right)  }=S_{1}S_{2}\ldots=S$ is the kneading sequence of
$+a$. Let $\mathtt{x}_{k}^{\left(  1\right)  }=\operatorname{It}\left(
x_{k}^{\left(  1\right)  }\right)  $ for $k=0,\ldots,p-1$. By symmetry
$\mathtt{x}_{1}^{\left(  1\right)  }=\tau S$ is the kneading sequence of $-a$.
It is also clear that $\mathtt{x}_{0}^{\left(  2\right)  }=\sigma^{p-1}\left(
S\right)  $ and $\mathtt{x}_{0}^{\left(  1\right)  }=\sigma^{p-1}\left(  \tau
S\right)  $. To each $k=1,\ldots,p$ corresponds a symbolic sequence
\[
\mathtt{x}_{k}^{\left(  2\right)  }=\sigma^{k-1}\left(  S\right)  .
\]
To each $k=1\ldots,p,$corresponds another sequence
\[
\mathtt{x}_{k}^{\left(  1\right)  }=\sigma^{k-1}\left(  \tau S\right)  .
\]
Naturally, we have
\[
\mathtt{x}_{k}^{\left(  1\right)  }=\tau\left(  \mathtt{x}_{k}^{\left(
2\right)  }\right)  .
\]
To each $z_{j}$ corresponds the symbolic itinerary $w_{j}=\operatorname{It}%
\left(  z_{j}\right)  $. We define $v_{j}=\operatorname{It}\left(
y_{j}\right)  $.

\begin{exmp}
\label{example1}Given the kneading sequence $\left(  RMB\right)  ^{\infty}$,
equivalent to $\left(  RMR\right)  ^{\infty}$ already used before. The
kneading pair is $\left(  \left(  RMB\right)  ^{\infty},\left(  LMA\right)
^{\infty}\right)  $, $\mathtt{x}_{1}^{\left(  2\right)  }=RMB=v_{2},$
$\mathtt{x}_{2}^{\left(  2\right)  }=MBR=v_{3}$, $\mathtt{x}_{0}^{\left(
2\right)  }=BRM=v_{1}$, $\mathtt{x}_{1}^{\left(  1\right)  }=LMA=v_{5},$
$\mathtt{x}_{2}^{\left(  1\right)  }=MAL=v_{6},$ and $\mathtt{x}_{0}^{\left(
1\right)  }=ALM=v_{4}$. We get $w_{1}=v_{5}$, $w_{2}=v_{4}$, $w_{3}=v_{3}$,
$w_{6}=v_{4}$, $w_{5}=v_{1}$ and $w_{6}=v_{2}$. The matrix $\pi$ is
\[
\pi=\left[
\begin{array}
[c]{rrrrrr}%
0 & 0 & 0 & 0 & 1 & 0\\
0 & 0 & 0 & 1 & 0 & 0\\
0 & 0 & 1 & 0 & 0 & 0\\
0 & 0 & 0 & 0 & 0 & 1\\
1 & 0 & 0 & 0 & 0 & 0\\
0 & 1 & 0 & 0 & 0 & 0
\end{array}
\right]  .
\]

\end{exmp}

\begin{exmp}
\label{example2}Given the kneading sequence $\left(  RLMB\right)  ^{\infty}$,
the kneading determinant $D\left(  t\right)  $ is such that
\begin{align*}
\left(  1+t\right)  D\left(  t\right)   &  =1+2\frac{-t-t^{2}+t^{4}}{1-t^{4}%
}\\
&  =\frac{1-2t-2t^{2}+t^{4}}{1-t^{4}}.
\end{align*}
The kneading pair is $\left(  \left(  RLMB\right)  ^{\infty},\left(
LRMA\right)  ^{\infty}\right)  $, $\mathtt{x}_{1}^{\left(  2\right)
}=RLMB=v_{2},$ $\mathtt{x}_{2}^{\left(  2\right)  }=LMBR=v_{3},$
$\mathtt{x}_{3}^{\left(  2\right)  }=MBRL=v_{4}$, $\mathtt{x}_{0}^{\left(
2\right)  }=BRLM=v_{1}$, $\mathtt{x}_{1}^{\left(  1\right)  }=LRMA=v_{6},$
$\mathtt{x}_{2}^{\left(  1\right)  }=RMAL=v_{7},$ $\mathtt{x}_{3}^{\left(
1\right)  }=MALR=v_{8}$ and $\mathtt{x}_{0}^{\left(  1\right)  }=ALRM=v_{5}$.
We get $w_{1}=v_{6}$, $w_{2}=v_{3}$, $w_{3}=v_{5}$, $w_{4}=v_{4}$,
$w_{5}=v_{8}$, $w_{6}=v_{1}$, $w_{7}=v_{7}$ and $w_{8}=v_{2}$. The Markov
matrix is
\[
\Psi=\left[
\begin{array}
[c]{rrrrrrr}%
0 & 0 & 0 & 1 & 1 & 1 & 0\\
1 & 1 & 1 & 0 & 0 & 0 & 0\\
0 & 0 & 0 & 0 & 0 & 1 & 1\\
0 & 0 & 1 & 1 & 1 & 0 & 0\\
1 & 1 & 0 & 0 & 0 & 0 & 0\\
0 & 0 & 0 & 0 & 1 & 1 & 1\\
0 & 1 & 1 & 1 & 0 & 0 & 0
\end{array}
\right]  .
\]
The matrix $\pi$ is
\[
\pi=\left[
\begin{array}
[c]{rrrrrrrr}%
0 & 0 & 0 & 0 & 0 & 1 & 0 & 0\\
0 & 0 & 1 & 0 & 0 & 0 & 0 & 0\\
0 & 0 & 0 & 0 & 1 & 0 & 0 & 0\\
0 & 0 & 0 & 1 & 0 & 0 & 0 & 0\\
0 & 0 & 0 & 0 & 0 & 0 & 0 & 1\\
1 & 0 & 0 & 0 & 0 & 0 & 0 & 0\\
0 & 0 & 0 & 0 & 0 & 0 & 1 & 0\\
0 & 1 & 0 & 0 & 0 & 0 & 0 & 0
\end{array}
\right]  .
\]

\end{exmp}

The matrix $\pi$ can also be used to reorder the shifts of the kneading
sequences, giving the vector $v=\left\{  v_{1},\ldots,v_{2p}\right\}  $ and
the vector $w=\left\{  w_{1},\ldots,w_{2p}\right\}  $, we have $w=\pi v$.

\section{Main Result}

\subsection{The Markov and kneading endomorphism in spaces of chain complexes}

Let $C_{0}$ be the vector space of the $0$-chains spanned by the shifts of the
kneading sequences $\left\{  v_{j}\right\}  _{j=1,\ldots,2p}$ this space is
isomorphic of the space of the $0$-chains spanned by the points of the orbit
$\left\{  y_{j}\right\}  _{j=1,\ldots,2p}$. The space $\pi\left(
C_{0}\right)  $ is spanned by $\left\{  w_{k}\right\}  _{k=1,\ldots,2p}$ which
is isomorphic to the space of the $0$-chains spanned by $\left\{
z_{j}\right\}  _{j=1,\ldots,2p}$. Let $C_{1}$ be the space of the $1$-chains
spanned by $\left\{  I_{k}\right\}  _{k=1,\ldots,2p-1}$, isomorphic to the
linear space of the $1$-chains spanned by $\left\{  I_{k}^{\prime}\right\}
_{k=1,\ldots,2p-1}$\ where $I_{k}^{\prime}$ is the set of all the admissible
sequences $w$: $w_{k}\preceq w\preceq w_{k+1}$. In what follows we identify
$I_{k}^{\prime}$ with $I_{k}$ and use the same symbol both for sequences and
intervals and call both the linear transformations and the corresponding
matrix representations by the same letters.

The border of a $1$-chain is obtained using the linear transformation
$\partial:C_{1}\rightarrow D_{0}$ such that $\partial\left(  I_{k}\right)
=w_{k+1}-w_{k}$, $\partial\left(  C_{1}\right)  =D_{0}$ where $D_{0}$ is
spanned by
\[
\left\{  w_{k+1}-w_{k}\right\}  _{k=1,\ldots,2p-1}.
\]
It is clear that $D_{0}\subset\pi\left(  C_{0}\right)  $.

We define the linear transformation $\partial_{s}:C_{1}\rightarrow D_{0}$ such
that
\begin{align*}
\partial_{s}\left(  I_{k}\right)   &  =\partial\left(  \mathbf{1}-\tau\right)
I_{k}\\
&  =\partial\left(  I_{k}-I_{2p-k}\right) \\
&  =\partial\left(  I_{k}\right)  -\partial\left(  I_{2p-k}\right)  .
\end{align*}
The image of $I_{k}$ by $\partial_{s}$ is
\[
\partial_{s}\left(  I_{k}\right)  =w_{k+1}-w_{k}-\left(  w_{2p+1-k}%
-w_{2p-k}\right)
\]
and is an element of $D_{0}$. We can define another linear transformation that
acts on $\pi\left(  C_{0}\right)  $ with matrix representation
\[
\mu=\left[
\begin{array}
[c]{rrrrrrr}%
-1 & 1 & 0 & \cdots & 0 & 1 & -1\\
0 & -1 & 1 & \cdots & 1 & -1 & 0\\
\vdots & \vdots & \vdots & \ddots & \vdots & \vdots & \vdots\\
0 & 1 & -1 & \cdots & -1 & 1 & 0\\
1 & -1 & 0 & \cdots & 0 & -1 & 1
\end{array}
\right]  .
\]
where $\mu_{ij}=\delta_{i+1,j}-\delta_{i,j}-\delta_{2p+1-i,j}+\delta
_{2p-i,j},$ $i=1,\ldots,2p-1$, $j=1,\ldots,2p$, and $\delta$ is the Kronecker
delta symbol. This linear transformation represents the order relation of the
points of the real line and the symmetry of the original mapping. It is
immediate from the above definitions that $Image\left(  \partial_{s}\right)
=B_{0}$ is a proper subspace of $D_{0}$ such that $D_{0}=B_{0}\oplus
\widetilde{B}$, where $\widetilde{B}$ has dimension one, and $B_{0}$ is
isomorphic to $\mu\pi\left(  C_{0}\right)  $. We define $\eta=\mu\pi$ and the
endomorphism$\ \omega$ acting on $C_{0}$ with matrix representation
\[
\omega=\left[
\begin{array}
[c]{cc}%
\sigma & 0\\
0 & \sigma
\end{array}
\right]  ,
\]
where $\sigma$ is the shift operator with matrix representation $p\times p$%
\[
\sigma=\left[
\begin{array}
[c]{rrrrr}%
0 & 1 & \cdots & 0 & 0\\
0 & 0 & \cdots & 0 & 0\\
\vdots & \vdots & \ddots & \vdots & \vdots\\
0 & 0 & \cdots & 0 & 1\\
1 & 0 & \cdots & 0 & 0
\end{array}
\right]  .
\]
Let $\alpha$ be the endomorphism induced in $B_{0}$ by the rotation $\omega$
in $C_{0}$ which results from the commutativity of the diagram
\[%
\begin{array}
[c]{ccccccc}
&  & _{\eta} &  & _{\partial_{s}} &  & \\
& C_{0} & \longrightarrow & B_{0} & \longleftarrow & C_{1} & \\
^{_{_{\omega}}} & \downarrow &  & \downarrow^{_{_{\alpha}}} &  & \downarrow &
^{_{_{\alpha}}}\\
& C_{0} & \longrightarrow & B_{0} & \longleftarrow & C_{1} & \\
&  & ^{\eta} &  & ^{\partial_{s}} &  &
\end{array}
.
\]
It is easy to see that $\alpha=-\Psi$, where $\Psi$ is the Markov matrix. Note
that $\eta\omega=\alpha\eta$. Every entry in the matrix $\alpha$ is
non-positive, because the images of the intervals are obtained by the images
of the boundary points and $F$ is reverse order in any of each interval of
continuity (lap).

\begin{exmp}
With the matrices of the examples \ref{example1} we have
\[
\eta=\left[
\begin{array}
[c]{rrrrrr}%
1 & -1 & 0 & 1 & -1 & 0\\
-1 & 0 & 1 & -1 & 0 & 1\\
0 & 0 & 0 & 0 & 0 & 0\\
1 & 0 & -1 & 1 & 0 & -1\\
-1 & 1 & 0 & -1 & 1 & 0
\end{array}
\right]
\]
and
\[
\eta\omega=\left[
\begin{array}
[c]{rrrrrr}%
0 & 1 & -1 & 0 & 1 & -1\\
1 & -1 & 0 & 1 & -1 & 0\\
0 & 0 & 0 & 0 & 0 & 0\\
-1 & 1 & 0 & -1 & 1 & 0\\
0 & -1 & 1 & 0 & -1 & 1
\end{array}
\right]  ,
\]
which is precisely $-\Psi\eta$.
\end{exmp}

\subsection{Matrix $\Theta$}

Giving the right $p$-periodic kneading sequence $S$, the symbolic itinerary of
the right discontinuity point is $\sigma^{p-1}S=S_{0}S_{1}S_{2}\ldots=S_{0}S$.
We construct the vector
\[
s\left(  S\right)  =\left[
\begin{array}
[c]{c}%
\Phi\left(  S_{0}\right) \\
\Phi\left(  S_{1}\right) \\
\vdots\\
\Phi\left(  S_{p-1}\right) \\
\Phi\left(  \tau S_{0}\right) \\
\Phi\left(  \tau S_{1}\right) \\
\vdots\\
\Phi\left(  \tau S_{p-1}\right)
\end{array}
\right]  =\left[
\begin{array}
[c]{c}%
1\\
\Phi\left(  S_{1}\right) \\
\vdots\\
\Phi\left(  S_{p-1}\right) \\
-1\\
\Phi\left(  \tau S_{1}\right) \\
\vdots\\
\Phi\left(  \tau S_{p-1}\right)
\end{array}
\right]  ,
\]
where $\Phi$ was defined in definition \ref{FI}. When applied to the other
kneading sequence $\tau S$ the vector takes the form
\[
s\left(  \tau S\right)  =\sigma^{p}s\left(  S\right)  .
\]
Let $\Gamma$ be a square matrix which columns $1$ and $p+1$ are $s$ and
$\sigma^{p}s$, respectively, and the other elements are zeros.

Now, we introduce the matrices $\gamma=\Gamma-I$ and $\Theta=\gamma\omega$.
The matrix $\Theta$ has the form
\[
^{\left[
\begin{array}
[c]{rrrrrrrrrrrr}%
^{0} & ^{0} & ^{0} & ^{\ldots} & ^{0} & ^{0} & ^{0} & ^{-\Phi\left(
S_{0}\right)  } & ^{0} & ^{\ldots} & ^{0} & ^{0}\\
^{0} & ^{\Phi\left(  S_{1}\right)  } & ^{-1} & ^{\ldots} & ^{0} & ^{0} & ^{0}
& ^{-\Phi\left(  S_{1}\right)  } & ^{0} & ^{\ldots} & ^{0} & ^{0}\\
^{0} & ^{\Phi\left(  S_{2}\right)  } & ^{0} & ^{\ldots} & ^{0} & ^{0} & ^{0} &
^{-\Phi\left(  S_{2}\right)  } & ^{0} & ^{\ldots} & ^{0} & ^{0}\\
^{\vdots} & ^{\vdots} & ^{\vdots} & ^{\ddots} & ^{\vdots} & ^{\vdots} &
^{\vdots} & ^{\vdots} & ^{\vdots} & ^{{}} & ^{\vdots} & ^{\vdots}\\
^{0} & ^{\Phi\left(  S_{p-2}\right)  } & ^{0} & ^{\ldots} & ^{0} & ^{-1} &
^{0} & ^{-\Phi\left(  S_{p-2}\right)  } & ^{0} & ^{\ldots} & ^{0} & ^{0}\\
^{-1} & ^{\Phi\left(  S_{p-1}\right)  } & ^{0} & ^{\ldots} & ^{0} & ^{0} &
^{0} & ^{-\Phi\left(  S_{p-1}\right)  } & ^{0} & ^{\ldots} & ^{0} & ^{0}\\
^{0} & ^{-\Phi\left(  S_{0}\right)  } & ^{0} & ^{\ldots} & ^{0} & ^{0} & ^{0}
& ^{0} & ^{0} & ^{\ldots} & ^{0} & ^{0}\\
^{0} & ^{-\Phi\left(  S_{1}\right)  } & ^{0} & ^{\ldots} & ^{0} & ^{0} & ^{0}
& ^{\Phi\left(  S_{1}\right)  } & ^{-1} & ^{\ldots} & ^{0} & ^{0}\\
^{0} & ^{-\Phi\left(  S_{2}\right)  } & ^{0} & ^{\ldots} & ^{0} & ^{0} & ^{0}
& ^{\Phi\left(  S_{2}\right)  } & ^{0} & ^{\ldots} & ^{0} & ^{0}\\
^{\vdots} & ^{\vdots} & ^{\vdots} & ^{{}} & ^{\vdots} & ^{\vdots} & ^{\vdots}
& ^{\vdots} & ^{\vdots} & ^{\ddots} & ^{\vdots} & ^{\vdots}\\
^{0} & ^{-\Phi\left(  S_{p-2}\right)  } & ^{0} & ^{\ldots} & ^{0} & ^{0} &
^{0} & ^{\Phi\left(  S_{p-2}\right)  } & ^{0} & ^{\ldots} & ^{0} & ^{-1}\\
^{0} & ^{-\Phi\left(  S_{p-1}\right)  } & ^{0} & ^{\ldots} & ^{0} & ^{0} &
^{-1} & ^{\Phi\left(  S_{p-1}\right)  } & ^{0} & ^{\ldots} & ^{0} & ^{0}%
\end{array}
\right]  .}%
\]

\begin{exmp}
With the same kneading sequences of the example \ref{example1}, we have
\[
\gamma=\left[
\begin{array}
[c]{rrrrrr}%
0 & 0 & 0 & -1 & 0 & 0\\
1 & -1 & 0 & -1 & 0 & 0\\
0 & 0 & -1 & 0 & 0 & 0\\
-1 & 0 & 0 & 0 & 0 & 0\\
-1 & 0 & 0 & 1 & -1 & 0\\
0 & 0 & 0 & 0 & 0 & -1
\end{array}
\right]
\]
and
\[
\Theta=\left[
\begin{array}
[c]{rrrrrr}%
0 & 0 & 0 & 0 & -1 & 0\\
0 & 1 & -1 & 0 & -1 & 0\\
-1 & 0 & 0 & 0 & 0 & 0\\
0 & -1 & 0 & 0 & 0 & 0\\
0 & -1 & 0 & 0 & 1 & -1\\
0 & 0 & 0 & -1 & 0 & 0
\end{array}
\right]  .
\]

\end{exmp}

\begin{prop}
The following diagram commutes
\[%
\begin{array}
[c]{ccccc}
&  & _{\eta} &  & \\
& C_{0} & \longrightarrow & B_{0} & \\
& ^{_{_{\gamma}}}\downarrow &  & \downarrow^{_{_{-I}}} & \\
& C_{0} & \longrightarrow & B_{0} & \\
&  & ^{\eta} &  &
\end{array}
.
\]

\end{prop}

\smallskip\noindent\textbf{Proof.} We must show that $\eta\gamma=-\eta$ or
$\eta\left(  \Gamma-I\right)  =-\eta$, in other words that $\eta\Gamma=0$, or
that $s\left(  S\right)  ,s\left(  \tau S\right)  \in kernel\left(
\eta\right)  $. But $\eta=\mu\pi$, and $\pi$ reorders $s\left(  S\right)  $ in
terms of the order of the real line, giving
\[
\pi s\left(  S\right)  =\left[
\begin{array}
[c]{l}%
\nu\left(  \mathbb{I}\left(  z_{1}\right)  \right) \\
\nu\left(  \mathbb{I}\left(  z_{2}\right)  \right) \\
\vdots\\
\nu\left(  \mathbb{I}\left(  z_{2p-1}\right)  \right) \\
\nu\left(  \mathbb{I}\left(  z_{2p}\right)  \right)
\end{array}
\right]  ,
\]
knowing that $\mathbb{I}\left(  z_{j}\right)  =\tau\mathbb{I}\left(
z_{2p-j}\right)  ,$ $j=1,\ldots p$, we have $\nu\left(  \mathbb{I}\left(
z_{j}\right)  \right)  =-\nu\left(  \mathbb{I}\left(  z_{2p-j}\right)
\right)  $. It is obvious that $\mu\pi s\left(  S\right)  =\mu\pi s\left(
\tau S\right)  =0$. $\square$

\begin{thm}
\label{MAIN}The characteristic polynomial of the matrix $\Theta=\gamma\omega$
is
\[
P_{\Theta}\left(  t\right)  =\det\left(  I-t\Theta\right)  =\left(  1-\left(
-1\right)  ^{p}t^{p}\right)  ^{2}\left(  1+t\right)  D\left(  t\right)  ,
\]
where $D\left(  t\right)  $ is the kneading determinant.
\end{thm}

\smallskip\noindent\textbf{Proof.} The determinant of the matrix $I-t\Theta$
is
\[
\left\vert
\begin{array}
[c]{cccccccccccc}%
1 & 0 & 0 & \ldots & 0 & 0 & 0 & t\Phi\left(  S_{0}\right)  & 0 & \ldots & 0 &
0\\
0 & 1-t\Phi\left(  S_{1}\right)  & t & \ldots & 0 & 0 & 0 & t\Phi\left(
S_{1}\right)  & 0 & \ldots & 0 & 0\\
0 & -t\Phi\left(  S_{2}\right)  & 1 & \ldots & 0 & 0 & 0 & t\Phi\left(
S_{2}\right)  & 0 & \ldots & 0 & 0\\
\vdots & \vdots & \vdots & \ddots & \vdots & \vdots & \vdots & \vdots & \vdots
&  & \vdots & \vdots\\
0 & -t\Phi\left(  S_{p-2}\right)  & 0 & \ldots & 1 & t & 0 & t\Phi\left(
S_{p-2}\right)  & 0 & \ldots & 0 & 0\\
t & -t\Phi\left(  S_{p-1}\right)  & 0 & \ldots & 0 & 1 & 0 & t\Phi\left(
S_{p-1}\right)  & 0 & \ldots & 0 & 0\\
0 & t\Phi\left(  S_{0}\right)  & 0 & \ldots & 0 & 0 & 1 & 0 & 0 & \ldots & 0 &
0\\
0 & t\Phi\left(  S_{1}\right)  & 0 & \ldots & 0 & 0 & 0 & 1-t\Phi\left(
S_{1}\right)  & t & \ldots & 0 & 0\\
0 & t\Phi\left(  S_{2}\right)  & 0 & \ldots & 0 & 0 & 0 & -t\Phi\left(
S_{2}\right)  & 1 & \ldots & 0 & 0\\
\vdots & \vdots & \vdots &  & \vdots & \vdots & \vdots & \vdots & \vdots &
\ddots & \vdots & \vdots\\
0 & t\Phi\left(  S_{p-2}\right)  & 0 & \ldots & 0 & 0 & 0 & -t\Phi\left(
S_{p-2}\right)  & 0 & \ldots & 1 & t\\
0 & t\Phi\left(  S_{p-1}\right)  & 0 & \ldots & 0 & 0 & t & -t\Phi\left(
S_{p-1}\right)  & 0 & \ldots & 0 & 1
\end{array}
\right\vert
\]
We first multiply the row $1$ by $-t$ and add the result to the row $p$, we do
the same with the rows $p+1$ and $2p$. Then we develop the determinant by the
columns $1$ and $p$ getting a $\left(  2p-1\right)  \times\left(  2p-1\right)
$ determinant:
\[
\left\vert
\begin{array}
[c]{cccccccccc}%
^{1-t\Phi\left(  S_{1}\right)  } & ^{t} & ^{\ldots} & ^{0} & ^{0} &
^{t\Phi\left(  S_{1}\right)  } & ^{0} & ^{\ldots} & ^{0} & ^{0}\\
^{-t\Phi\left(  S_{2}\right)  } & ^{1} & ^{\ldots} & ^{0} & ^{0} &
^{t\Phi\left(  S_{2}\right)  } & ^{0} & ^{\ldots} & ^{0} & ^{0}\\
^{\vdots} & ^{\vdots} & ^{\ddots} & ^{\vdots} & ^{\vdots} & ^{\vdots} &
^{\vdots} &  & ^{\vdots} & ^{\vdots}\\
^{-t\Phi\left(  S_{p-2}\right)  } & ^{0} & ^{\ldots} & ^{1} & ^{t} &
^{t\Phi\left(  S_{p-2}\right)  } & ^{0} & ^{\ldots} & ^{0} & ^{0}\\
^{-t\Phi\left(  S_{p-1}\right)  } & ^{0} & ^{\ldots} & ^{0} & ^{1} &
^{t\Phi\left(  S_{p-1}\right)  -t^{2}\Phi\left(  S_{0}\right)  } & ^{0} &
^{\ldots} & ^{0} & ^{0}\\
^{t\Phi\left(  S_{1}\right)  } & ^{0} & ^{\ldots} & ^{0} & ^{0} &
^{1-t\Phi\left(  S_{1}\right)  } & ^{t} & ^{\ldots} & ^{0} & ^{0}\\
^{t\Phi\left(  S_{2}\right)  } & ^{0} & ^{\ldots} & ^{0} & ^{0} &
^{-t\Phi\left(  S_{2}\right)  } & ^{1} & ^{\ldots} & ^{0} & ^{0}\\
^{\vdots} & ^{\vdots} &  & ^{\vdots} & ^{\vdots} & ^{\vdots} & ^{\vdots} &
^{\ddots} & ^{\vdots} & ^{\vdots}\\
^{t\Phi\left(  S_{p-2}\right)  } & ^{0} & ^{\ldots} & ^{0} & ^{0} &
^{-t\Phi\left(  S_{p-2}\right)  } & ^{0} & ^{\ldots} & ^{1} & ^{t}\\
^{t\Phi\left(  S_{p-1}\right)  -t^{2}\Phi\left(  S_{0}\right)  } & ^{0} &
^{\ldots} & ^{0} & ^{0} & ^{-t\Phi\left(  S_{p-1}\right)  } & ^{0} & ^{\ldots}
& ^{0} & ^{1}%
\end{array}
\right\vert ,
\]
then we multiply the row $p-1$ by $-t$ and add the result to the row $p-2$. We
do the same with the last two rows. Then we develop the determinant by the
columns $p-1$ and the last one getting a $\left(  2p-2\right)  \times\left(
2p-2\right)  $ determinant:
\[
\left\vert
\begin{array}
[c]{cccccccc}%
^{1-t\Phi\left(  S_{1}\right)  } & ^{t} & ^{\ldots} & ^{0} & ^{t\Phi\left(
S_{1}\right)  } & ^{0} & ^{\ldots} & ^{0}\\
^{-t\Phi\left(  S_{2}\right)  } & ^{1} & ^{\ldots} & ^{0} & ^{t\Phi\left(
S_{2}\right)  } & ^{0} & ^{\ldots} & ^{0}\\
^{\vdots} & ^{\vdots} & ^{\ddots} & ^{\vdots} & ^{\vdots} & ^{\vdots} & ^{{}}
& ^{\vdots}\\
^{r_{1}\left(  t\right)  } & ^{0} & ^{\ldots} & ^{1} & ^{r_{2}\left(
t\right)  } & ^{0} & ^{\ldots} & ^{0}\\
^{t\Phi\left(  S_{1}\right)  } & ^{0} & ^{\ldots} & ^{0} & ^{1-t\Phi\left(
S_{1}\right)  } & ^{t} & ^{\ldots} & ^{0}\\
^{t\Phi\left(  S_{2}\right)  } & ^{0} & ^{\ldots} & ^{0} & ^{-t\Phi\left(
S_{2}\right)  } & ^{1} & ^{\ldots} & ^{0}\\
^{\vdots} & ^{\vdots} & ^{{}} & ^{\vdots} & ^{\vdots} & ^{\vdots} & ^{\ddots}
& ^{\vdots}\\
^{r_{2}\left(  t\right)  } & ^{0} & ^{\ldots} & ^{0} & ^{r_{1}\left(
t\right)  } & ^{0} & ^{\ldots} & ^{1}%
\end{array}
\right\vert \text{,}%
\]
where $r_{1}\left(  t\right)  =-t\Phi\left(  S_{p-2}\right)  +t^{2}\Phi\left(
S_{p-1}\right)  $ and $r_{2}\left(  t\right)  =t\Phi\left(  S_{p-2}\right)
-t^{2}\Phi\left(  S_{p-1}\right)  +t^{3}\Phi\left(  S_{0}\right)  $. Repeating
this reducing process we get a $2\times2$ determinant
\[
\left\vert
\begin{array}
[c]{ll}%
1-\sum_{k=1}^{p-1}\left(  -1\right)  ^{k}t^{k}\Phi\left(  S_{k}\right)  &
-\sum_{k=1}^{p}\left(  -1\right)  ^{k}t^{k}\Phi\left(  S_{k}\right) \\
-\sum_{k=1}^{p}\left(  -1\right)  ^{k}t^{k}\Phi\left(  S_{k}\right)  &
1-\sum_{k=1}^{p-1}\left(  -1\right)  ^{k}t^{k}\Phi\left(  S_{k}\right)
\end{array}
\right\vert .
\]
Remembering that
\[
u_{p}\left(  t\right)  =\sum_{k=1}^{p}\left(  -1\right)  ^{k}t^{k}\Phi\left(
S_{k}\right)  \text{,}%
\]%
\[
S_{0}=S_{p}=B\text{ and }\left(  -1\right)  ^{p}t^{p}\Phi\left(  B\right)
=\left(  -1\right)  ^{p}t^{p}\text{,}%
\]
the previous determinant is equal to
\[
\left\vert
\begin{array}
[c]{cc}%
\left(  1-\left(  -1\right)  ^{p}t^{p}\right)  +u_{p}\left(  t\right)  &
-u_{p}\left(  t\right) \\
-u_{p}\left(  t\right)  & \left(  1-\left(  -1\right)  ^{p}t^{p}\right)
+u_{p}\left(  t\right)
\end{array}
\right\vert ,
\]
which gives
\[
P_{\Theta}\left(  t\right)  =\left(  1-\left(  -1\right)  ^{p}t^{p}\right)
^{2}\left(  1+2\frac{u_{p}\left(  t\right)  }{1-\left(  -1\right)  ^{p}t^{p}%
}\right)
\]
and this is precisely
\[
P_{\Theta}\left(  t\right)  =\left(  1-\left(  -1\right)  ^{p}t^{p}\right)
^{2}\left(  1+t\right)  D\left(  t\right)  ,
\]
as desired. $\square$

\begin{exmp}
We use the kneading sequences of the example \ref{example2} to illustrate this
last result, the matrix $\Theta$ is
\[
\left[
\begin{array}
[c]{rrrrrrrr}%
0 & 0 & 0 & 0 & 0 & -1 & 0 & 0\\
0 & 1 & -1 & 0 & 0 & -1 & 0 & 0\\
0 & -1 & 0 & -1 & 0 & 1 & 0 & 0\\
-1 & 0 & 0 & 0 & 0 & 0 & 0 & 0\\
0 & -1 & 0 & 0 & 0 & 0 & 0 & 0\\
0 & -1 & 0 & 0 & 0 & 1 & -1 & 0\\
0 & 1 & 0 & 0 & 0 & -1 & 0 & -1\\
0 & 0 & 0 & 0 & -1 & 0 & 0 & 0
\end{array}
\right]  ,
\]
with characteristic polynomial
\[
\left(  1-t^{4}\right)  \left(  1-2t-2t^{2}+t^{4}\right)  ,
\]
which agrees with the value of the kneading determinant $\left(  1+t\right)
D\left(  t\right)  =\frac{1-2t-2t^{2}+t^{4}}{1-t^{4}}$.
\end{exmp}

We have now all the ingredients to state the main results of this paper.

\begin{thm}
The following diagram commutes
\[%
\begin{array}
[c]{ccccc}
&  & _{\eta} &  & \\
& C_{0} & \longrightarrow & B_{0} & \\
& \downarrow^{_{_{\Theta}}} &  & \downarrow^{_{_{\Psi}}} & \\
& C_{0} & \longrightarrow & B_{0} & \\
&  & ^{\eta} &  &
\end{array}
.
\]
and $P_{\Theta}\left(  t\right)  =\left(  1+t\right)  \det\left(
I-t\Psi\right)  $.
\end{thm}

\smallskip\noindent\textbf{Proof.} Noticing that $\Theta=\gamma\omega$ and
$\Psi=-\alpha$ the result is only a direct consequence of
\[%
\begin{array}
[c]{ccc}
& _{_{\eta}} & \\
C_{0} & \longrightarrow & B_{0}\\
\downarrow_{^{\omega}} & _{_{\eta}} & \downarrow_{^{\alpha}}\\
C_{0} & \longrightarrow & B_{0}\\
\downarrow^{_{_{\gamma}}} & _{_{\eta}} & \downarrow^{_{_{-I}}}\\
C_{0} & \longrightarrow & B_{0}%
\end{array}
\]
conjugated with the fact that the two rows in the next diagram
\[%
\begin{array}
[c]{ccccccccc}
&  &  & _{\operatorname{inj}} &  & _{_{\eta}} &  &  & \\
0 & \longrightarrow & \widetilde{B} & \longrightarrow & C_{0} &
\longrightarrow & B_{0} & \longrightarrow & 0\\
&  & \downarrow^{_{_{-I}}} & _{\operatorname{inj}} & \downarrow^{_{_{\Theta}}}
& _{_{\eta}} & \downarrow^{_{_{\Psi}}} &  & \\
0 & \longrightarrow & \widetilde{B} & \longrightarrow & C_{0} &
\longrightarrow & B_{0} & \longrightarrow & 0
\end{array}
\]
are exact sequences, where $\operatorname{inj}$ is the natural\ embedding.
$\square$

\begin{cor}
The inverse of the least root in the unit interval of the periodic kneading
determinant is the spectral radius of the Markov matrix.
\end{cor}

\smallskip\noindent\textbf{Proof.} Is an immediate consequence of the
relation
\[
\left(  1-\left(  -1\right)  ^{p}t^{p}\right)  ^{2}D\left(  t\right)
=\det\left(  I-t\Psi\right)  ,
\]
obtained in the last theorem. $\square$

We think that the results of this work can be extended to general
discontinuous maps with finite number of discontinuities. That is a natural
extension of this work.

\section*{Acknowledgement}

The author thanks some precious comments of the referee that simplified and
improved a great deal the scope and clarity of this work. Namely, the
possibility of infinite jumps at the discontinuity points. The author was
partially funded by FCT/Portugal through project PEst-OE/EEI/LA0009/2013.

\emergencystretch=\hsize

\begin{center}
\rule{6 cm}{0.02 cm}
\end{center}


\begin{thebibliography}{99}                                                                                               %


\bibitem {Al}L.~Alsed\'a and F.~Ma{\~n}osas. \newblock Kneading theory and
rotation intervals for a class of circle maps of degree one. \newblock {\em
Nonlinearity}, 3(2):413--452, 1990.

\bibitem {Alv}J.~F. Alves and J.~Sousa-Ramos. \newblock Kneading theory for
tree maps. \newblock {\em Ergodic Theory and Dynamical Systems},
24(4):957--985, 2004.

\bibitem {CUR}J.~H. Curry, L.~Garnett, and D.~Sullivan. \newblock On the
iteration of a rational function: computer experiments with newton's method.
\newblock {\em Communications in mathematical physics}, 91(2):267--277, 1983.

\bibitem {EH}E.~D'Aniello and H.~M. Oliveira. \newblock Pitchfork bifurcation
for non-autonomous interval maps. \newblock {\em Difference
Equations and Applications}, 15(3):291--302, 2009.

\bibitem {dev-keen}R.~Devaney and L.~Keen. \newblock {\em Dynamics of
Tangent}, volume 1342 of \emph{Lecture Notes in Mathematics}, pages 105--111.
\newblock Springer, Berlin, New York, 1988.

\bibitem {DS}W.~de~Melo and S.~Strien. \newblock {\em One-dimensional
dynamics}. \newblock Springer, Berlin, Heildelberg, 1993.

\bibitem {Glen}P.~Glendinning and C.~Sparrow. \newblock Prime and
renormalisable kneading invariants and the dynamics of expanding lorenz maps.
\newblock {\em Physica D: Nonlinear Phenomena}, 62(1):22--50, 1993.

\bibitem {keen-kotus}L.~Keen and J.~Kotus. \newblock Dynamics of the family
$\lambda tan z$. \newblock {\em Conformal Geometry and Dynamics}, 1(1):28--57, 1997.

\bibitem {LSR}J.~P. Lampreia and J.~Sousa-Ramos. \newblock Symbolic dynamics
of bimodal maps. \newblock {\em Portugaliae Mathematica}, 54(1):1--18, 2013.

\bibitem {knead}J.~Milnor and W.~Thurston. \newblock {\em On iterated maps
of the interval}, volume 1342 of \emph{Lecture Notes in Mathematics}, pages
465--563. \newblock Springer, Berlin, 1988.

\bibitem {OLI}H.~Oliveira and J.~Sousa-Ramos. \newblock Iterates of
transcendent meromorphic maps. \newblock {\em Grazer Mathematische Berichte},
(346):313--321, 2004.
\end{thebibliography}
\end{document}